\documentclass[twoside]{article}

\usepackage{amsmath,amsthm,amssymb}
\usepackage[colorlinks=false,linkcolor=blue,citecolor=blue,urlcolor=blue]{hyperref}
\usepackage{fancynum}\setfnumgsym{\hskip 2.2pt}

\theoremstyle{theorem}
\newtheorem{theorem}{Theorem}[section]
\newtheorem*{thm*}{Theorem}

\numberwithin{equation}{section}
\theoremstyle{definition}
\newtheorem{remark}[theorem]{Remark}
\newtheorem{example}[theorem]{Example}

\newcommand{\cp}{\mathbb{C}P}
\newcommand{\zz}{\mathbb{Z}}

\title{\Large \bf Geometrical realization of low dimensional complete intersections$^{\ast}$}
\author{\large \bf Jianbo Wang$^{1}$, Jianpeng Du$^{2}$}
\date{}

\renewenvironment{abstract}{
        \small
        \quotation
         \noindent {\bfseries \abstractname } }
      {\if@twocolumn\else\endquotation\fi}

\def\Subsec{\@StartSubsection{subsection}{2}{\z@}
                                     {-3.25ex\@plus -1ex \@minus -.2ex}
                                     {1.5ex \@plus .2ex}
                                     {\normalfont\normalsize\bfseries\boldmath}}
                                     
\begin{document}

\maketitle

\renewcommand{\thefootnote}{\fnsymbol{footnote}}

\footnotetext{\hskip -5.5mm
$^1$  Department of Mathematics, School of Science, Tianjin University, Weijin Road 92, Nankai District, Tianjin  300072, P.R.China. {E-mail:} wjianbo@tju.edu.cn\\
$^{2}$  Graduate School, China Academy of Engineering Physics, P.O.Box 2101, Garden Road 6, Haidian District, Beijing 100088, P.R.China. {E-mail:} 1261708156@qq.com\\
 $^{\ast}$  The first author is supported by NSFC grant No.11001195 and Beiyang Elite Scholar Program of Tianjin University, No. 60301016.}

\begin{abstract}
This paper aims to give some examples of diffeomorphic (or homeomorphic) low dimensional complete intersections, which can be considered as a geometrical realization of classification theorems about complete intersections. A conjecture of Libgober and Wood \cite[$\S 9$]{LW82}(Topology. {\bf 21}, 1982, 469--482)
will be confirmed by one of examples.


 \noindent{\bf 2000 MR Subject Classification}
14M10, 57R19
\end{abstract}

\renewcommand{\thefootnote}{\arabic{footnote}}
\baselineskip 14pt

\section{Introduction}

Let $X_n(\underline{d})\subset \cp^{n+r}$ be a smooth complete
intersection of  multidegree $\underline{d}:=(d_1,\dots,d_r)$, the product $d_1d_2\cdots d_r$ is called
total degree, denoted by $d$. The celebrated Lefschetz hyperplane section theorem asserts that the
pair $(\cp^{n+r},X_n(\underline{d}))$ is $n$-connected. Thom
showed that the complex dimension $n$ and the multidegree $\underline{d}$
determine the diffeomorphism type of $X_n(\underline{d})$. However, the
multidegree is not a topological invariant. An interesting and
challenging problem is to classify complete intersections by
numerical topological invariants, such as total degree, the Pontrjagin classes and Euler characteristics, which usually are polynomials on
the multidegree (see Section \ref{cha}). 

In dimension $1$, the classification of complete intersections follows from the classical
theory of Riemann surfaces. In dimension $3$, appealing to general classification theorems in
differential topology, the classification is
established by Jupp \cite{J} and Wall \cite{Wa}, while in dimension $2$
the homeomorphism classification has been settled by Freedman \cite{Fr}.
Ebeling \cite{E} and Libgober-Wood \cite{LW90} independently found examples of homeomorphic complex $2$-dimensional complete intersections but not diffeomorphic.
In \cite{FK}, Fang and Klaus  proved that, in dimensions $n=4$,
two complete intersections $X_n(\underline{d})$ and $X_n(\underline{d}^\prime)$ are homeomorphic if and only if they have the same total degree, Pontrjagin classes and Euler characteristics.
Furthermore, Fang and the first author \cite{FW} generalized this homeomorphism result to dimension $n=5,6,7$.

\begin{theorem}[Fang-Klaus $\&$ Fang-Wang]\label{FKW} Two complete intersections $X_n(\underline{d})$ and
$X_n(\underline{d}^\prime)$ are homeomorphic if and only if they
have the same total degree, Pontrjagin classes and Euler
characteristics, provided $n=4,5,6,7$.
\end{theorem}

With the help of Kreck's modified surgery theory \cite{Kr}, Traving \cite{Tr} obtained
partial classification results in higher dimensions under some
restrictions on the total degree. Particularly, to the prime factorization of total degree $d=\prod_{p~\textrm{primes}}p^{\nu_{p}(d)}$, under the assumption that
$\nu_p(d)\geqslant \frac{2n+1}{2(p-1)}+1$ for all primes $p$ with $p(p-1)\leqslant n+1$, Traving proved the following result (see \cite[Theorem A]{Kr} or \cite{Tr}).
\begin{theorem}[Traving]\label{Traving}
Two complete intersections $X_n(\underline{d})$ and $X_n(\underline{d}^\prime)$
of complex dimension $n>2$ fulfilling the assumption above for the total degree are
diffeomorphic if and only if the total degrees, the Pontrjagin classes and the
Euler characteristics agree.
\end{theorem}

It is well known all complete intersections of fixed multidegree are diffeomorphic. On the other hand, there exist diffeomorphic complete intersections with different multidegrees. For lower dimensions, such as complex dimensions $2,3,4,5$, the diffeomorphic examples can be found in \cite{B,B00,LW82,Wang}.

The aim of this paper is to give examples of diffeomorphic (or homeomorphic) complex $n$-dimensional complete intersections,
mainly $n\leqslant 7$, which can be considered as a geometrical realization of above Theorems \ref{FKW} and \ref{Traving}. On the other hand, Libgober and Wood \cite[$\S 9$]{LW82} conjectured that  there exist the diffeomorphic complex 3-dimensional complete intersections with different first Chern classes, which give an example of a disconnected moduli space. In this paper, the above  conjecture will be confirmed by one of the listed examples.

These examples are partially found by computer searching. Note that, according to \cite[Proposition 7.3]{LW82}, if $X_n(d_1,\dots,d_r)\subset \cp^{n+r}$ satisfies $n>2$ and $r\leqslant \frac{n+2}{2}$, then the total degree and Pontrjagin classes of $X_n(d_1,\dots,d_r) $ determine the multidegree.  Thus, it is impossible to find out such an homeomorphic or diffeomorphic example in which the codimension $r$ is small relative to complex dimension $n$. It is worth stating that the degrees and codimensions in our examples will be as smaller as possible.

\section{Characteristic classes of complete intersections}\label{cha}

For a complete intersection
$X_n(\underline{d})$, let $H$ be the restriction of the dual bundle of the canonical line bundle over $\cp^{n+r}$ to $X_n(\underline{d})$, and $x= c_1(H)\in H^2(X_n(\underline{d});\zz)$.
Associate the multidegree $\underline{d}=(d_1,d_2,\dots,d_r)$, define the power sums $s_i=\sum_{j=1}^{r}d_j^i$ for $1\leqslant i\leqslant n$.
Then the Chern classes, Pontrjagin classes and Euler characteristic are presented as follows (\cite{LW82}):
\begin{align*}
c_k & =\frac{1}{k!}g_k(n+r+1-s_1,\dots,n+r+1-s_k)x^k, 1\leqslant k\leqslant n;\\
p_k & =\frac{1}{k!}g_k(n+r+1-s_2,\dots,n+r+1-s_{2k})x^{2k}, 1\leqslant k\leqslant
\left[\frac{n}{2}\right];\\
e & =c_n(X_n(\underline{d}))\cap [X_n(\underline{d})]
=d\frac{1}{n!}g_n(n+r+1-s_1,\dots,n+r+1-s_n).
\end{align*}
Where $g_k (k\geqslant 1)$ are the polynomials that can be iteratively computed from the Newton formula:
\begin{equation*}
s_k-g_1(s_1)s_{k-1}+\frac{1}{2}g_2(s_1,s_2)s_{k-2}+\cdots+(-1)^k\frac{1}{k!}g_k(s_1,s_2,\dots,s_k) k=0.
\end{equation*}
For example, here are the first seven $g_k$'s.
\begin{align*}
g_1(s_1) & =s_1,\\
g_2(s_1,s_2) & =s_1^2-s_2,\\
g_3(s_1,s_2,s_3) & =s_1^3-3s_1s_2+2s_3,\\
g_4(s_1,\dots,s_4) & =s_1^4-6s_1^2s_2+8s_1s_3+3s_2^2-6s_4,\\
g_5(s_1,\dots,s_5) & =s_1^5-10s_1^3s_2+20s_1^2s_3-30s_1s_4+15s_1s_2^2-20s_2s_3+24s_5,\\
g_6(s_1,\dots,s_6) & =s_1^6-15s_1^4s_2+40s_1^3s_3-90s_1^2s_4+45s_1^2s_2^2-120s_1s_2s_3\\
& \quad +144s_1s_5-15s_2^3+90s_2s_4+40s_3^2-120s_6.\\
g_7(s_1,\dots,s_7) & =s_1^7-21s_1^5s_2+70s_1^4s_3-210s_1^3s_4+105s_1^3s_2^2-420s_1^2s_2s_3\\
& \quad +504s_1^2s_5-105s_1s_2^3+630s_1s_2s_4+280s_1s_3^2-840s_1s_6\\
& \quad +210s_2^2s_3-504s_2s_5-420s_3s_4+720s_7.
\end{align*}

Note that the $k^{\textrm{th}}$ Chern classes $c_k$ and Pontrjagin class $p_k$ are integral multiples of $x^k$ and $x^{2k}$, so we can compare this invariant for different
complete intersections. For convenience, throughout the rest of the paper,
the Chern class $c_k$ and Pontrjagin class $p_k$ of $X_n(\underline{d})$ are viewed as the multiple of $x^k$ and $x^{2k}$ respectively.

\section{Examples of diffeomorphic (homeomorphic) low dimensional complete intersections}

In this section, some examples of diffeomorphic or homeomorphic complex $n$-dimensional complete intersections will be listed for $n=2,3,\dots,7$, and one of examples will confirm the conjecture of Libgober and Wood \cite[$\S 9$]{LW82}.

\subsection{Complex $2$-dimensional complete intersections}
For complex $2$-dimensional complete intersection $X_2(d_1,\dots,d_r)$, the total degree, first Chern class,
Pontrjagin class and Euler characteristic are as follows:
\begin{align*}
 d & =d_1\times\cdots\times  d_r, \\
 c_1 & =3+r-s_1, \\
 p_1 & =3+r-s_2, \\
 e & = \frac{d}{2}\left[(3+r-s_1)^2-(3+r-s_2)\right].
\end{align*}
Two complete intersections $X_2(d_1,\dots,d_r)$ and $X_2(d_1^\prime,\dots,d_s^\prime)$, which satisfy the following conditions,
\begin{align*}
 d\cdot p_1 & = d^\prime\cdot p_1^\prime, \\
 e & = e^\prime,\\
 c_1 & \equiv c_1^\prime \pmod 2, \\
 c_1 & \ne c_1^\prime.
\end{align*}
will give an example of homeomorhpic but non-diffeomorphic complex $2$-dimensional complete intersections (see \cite{E,LW90}).
Note that,  Ebeling \cite{E} and Libgober-Wood \cite{LW90} independently found that
$X_2(10,7,7,6,6,3,3)$ and $X_2(9,5,3,3,3,3,3,2,2)$ are
homeomorphic but not diffeomorphic. Here, some other examples with lower codimensions are listed.

\begin{example}The following invariants give seven examples of pairs of complex 2-dimensional complete intersections which are homeomorphic but not diffeomorphic. The two smoothings in each pair have distinct $c_1$'s.
\begin{table}[h]
\centerline{\(
\begin{array}{l|lll} \hline
X_2(\underline{d}) & d\cdot p_1  & e & c_1\\ \hline
X_2(6,5,3) & -5760  & 5760 & -8\\
X_2(5,2,2,2,2,2) & -5760 & 5760 & -6\\  \hline
X_2(14,14,8) & -705600  & 1058400 & -30\\
X_2(14,7,5,5) & -705600 & 1058400 & -24\\  \hline
X_2(15,11,10,3) & -2217600  & 3643200 & -32\\
X_2(11,10,5,2,2,2,2) & -2217600 & 3643200 & -24\\  \hline
X_2(10,9,7,7,7) & -9878400 & 20744640 & -32\\
X_2(7,7,7,5,2,2,2,2,2) & -9878400 & 20744640 & -24\\ \hline
X_2(10,9,9,6,2) & -2857680 & 5239080 & -28\\
X_2(10,7,7,3,3,3) & -2857680 & 5239080 & -24\\ \hline
X_2(15,13,7,3,2) & -3669120 & 6027840 & -32\\
X_2(13,7,5,2,2,2,2,2) & -3669120 & 6027840 & -24\\ \hline
X_2(10,7,7,6,3,3) & -6429780 & 12859560 & -27\\
X_2(9,5,3,3,3,3,3,2,2) & -6429780 & 12859560 & -21\\ \hline
\end{array}
\)}
\vskip .2cm
\caption{Homeomorphic but non-diffeomorphic 2-dim complete intersections}
\end{table}
\end{example}

\subsection{Complex $3$-dimensional complete intersections}

For complex $3$-dimensional complete intersection $X_3(d_1,\dots,d_r)$, the related topological characteristic classes are as follows:
\begin{align*}
 d & =d_1\times\cdots\times  d_r, \\
 p_1 & =4+r-s_2, \\
 e & = \frac{d}{6}\left[(4+r-s_1)^3-3(4+r-s_1)(4+r-s_2)+2(4+r-s_3)\right],\\
 c_1 & =4+r-s_1.
\end{align*}
Tow complex $3$-dimensional complete intersections are diffeomorphic, if and only if they have the same total degree $d$, the first Pontrjagin class $p_1$ and
Euler characteristics(\cite{J,Wa}). In \cite[$\S 9$]{LW82}, Libgober and Wood conjectured that there exist such examples of diffeomorphic complex 3-dimensional complete intersections with different first Chern classes, which would give an example of a disconnected moduli space. Br{\"u}ckmann \cite{B} show that there exist the diffeomorphic complex 3-dimensional complete intersections belonging to components of the moduli space of different dimensions, but with the same first Chern classes. In the following, some examples of diffeomorphic complex 3-dimensional complete intersections will be listed, and the conjecture in \cite{LW82} is confirmed by Table \ref{Diff3nonc1}.
\begin{example}
The following invariants give some pairs of diffeomorphic complex 3-dimensional complete intersections, where the ones with different first Chern classes $c_1$ give an example of a disconnected moduli space.
\begin{table}[h]
\centerline{\(
\begin{array}{l|llll} \hline
X_3(\underline{d}) & d & p_1  & e & c_1\\ \hline
X_3(20,20,11,7,4) & 123200 & -977 & -6974721600 & -53\\
X_3(22,16,14,5,5) & 123200 & -977 & -6974721600 & -53\\  \hline
X_3(14,14,5,4,4,4) & 62720 & -455 & -1068748800 & -35\\
X_3(16,10,7,7,2,2,2) & 62720 & -455 & -1068748800 & -35\\  \hline
\end{array}
\)}
\vskip .2cm
\caption{Diffeomorphic 3-dim complete intersections with same $c_1$}
\end{table}

\begin{table}[h]
\centerline{\(
\begin{array}{l|llll} \hline
X_3(\underline{d}) & d & p_1  & e & c_1\\ \hline
X_3(70,16,16,14,7,6) & 10536960 & -5683 & -7767425433600 & -119\\
X_3(56,49,8,6,5,4,4) & 10536960 & -5683 & -7767425433600 & -121\\  \hline
X_3(88,28,19,14,6,6) & 23595264 & -9147 & -35445749391360 & -151\\
X_3(76,56,11,7,6,6,2) & 23595264 & -9147 & -35445749391360 & -153\\  \hline
X_3(84,29,25,25,18,7) & 191835000 & -9510 & -384536710530000 & -178\\
X_3(60,58,49,9,5,5,5) & 191835000 & -9510 & -384536710530000 & -180\\  \hline
\end{array}
\)}
\vskip .2cm
\caption{Diffeomorphic 3-dim complete intersections with different $c_1$}\label{Diff3nonc1}
\end{table}
\end{example}

\subsection{Complex $4$-dimensional complete intersections}

For complex $4$-dimensional complete intersection $X_4(d_1,\dots,d_r)$, we have the following invariants:
\begin{align*}
 d & =d_1\times\cdots\times  d_r, \\
 p_1 & =5+r-s_2, \\
 p_2 & =\frac{1}{2}\left[(5+r-s_2)^2-(5+r-s_4)\right],  \\
 e & = \frac{d}{24}\big[(5+r-s_1)^4-6(5+r-s_1)^2(5+r-s_2)+8(5+r-s_1)(5+r-s_3)\\
& \quad +3(5+r-s_2)^2-6(5+r-s_4)\big]
\end{align*}
According to Theorem \ref{FKW}, to construct homeomorphic complex $4$-dimensional
complete intersections, we need to find out two different multidegrees to satisfy the above invariants all agree. Furthermore, by Theorem \ref{Traving}, if $\nu_2(d)>5.5$, then the
homeomorphic 4-dimensional complete intersections are diffeomorphic.

\begin{example}
Table \ref{Homeo4} give two pairs of homeomorphic complex 4-dimensional complete intersections.
\begin{table}[h]
\centerline{\(
\begin{array}{l|llll} \hline
X_4(\underline{d}) & d & p_1 & p_2 & e \\ \hline
X_4(66,63,29,23,6,4) & 66561264 & -9736 & 65253028 & 11837353833553248\\
X_4(69,58,36,14,11,3) & 66561264 & -9736 & 65253028 & 11837353833553248\\  \hline
X_4(46,44,33,27,17,15,10) & 4598629200 & -6472 & 25986916 & 546159737882484000\\
X_4(45,45,34,23,22,12,11) & 4598629200 & -6472 & 25986916 & 546159737882484000\\  \hline
\end{array}
\)}
\vskip .2cm
\caption{Homeomorphic 4-dim complete intersections}\label{Homeo4}
\end{table}
\end{example}

\begin{example}
 Table \ref{Diff4} give two pairs of diffeomorphic complex 4-dimensional complete intersections, since the total degree satisfy
\begin{align*}
& 488980800=2^6\times 3^4\times 5^2\times 7^3\times 11, \\
& 3953664000=2^{13} \times 3^3\times 5^3\times 11 \times 13.
\end{align*}
\begin{table}[h]
\centerline{\(
\begin{array}{l|llll} \hline
X_4(\underline{d}) & d & p_1 & p_2 & e \\ \hline
X_4(36,33,30,20,14,7,7) & 488980800 & -3967 & 9807916 & 19704249035856000\\
X_4(35,35,28,22,12,9,6) & 488980800 & -3967 & 9807916 & 19704249035856000\\  \hline
X_4(52,44,36,25,20,12,8) & 3953664000 & -7157 & 32268711 & 546278189783040000\\
X_4(50,48,32,26,22,10,9) & 3953664000 & -7157 & 32268711 & 546278189783040000\\  \hline
\end{array}
\)}
\vskip .2cm
\caption{Diffeomorphic 4-dim complete intersections}\label{Diff4}
\end{table}
\end{example}

\subsection{Complex $5$-dimensional complete intersections}

For complex $5$-dimensional complete intersection $X_5(d_1,\dots,d_r)$, we have the following invariants:
\begin{align*}
 d & =d_1\times\cdots\times  d_r, \\
 p_1 & =6+r-s_2, \label{p1} \\
 p_2 & =\frac{1}{2}{\big [}(6+r-s_2)^2-(6+r-s_4)\big{]}, \\
 e & =\frac{1}{5!}d\big{[}(6+r-s_1)^5-10(6+r-s_1)^3(6+r-s_2)+20(6+r-s_1)^2(6+r-s_3) \\
 & \quad -30(6+r-s_1)(6+r-s_4)
 +15(6+r-s_1)(6+r-s_2)^2  \\
  & \quad -20(6+r-s_2)(6+r-s_3)+24(6+r-s_5)\big ],
\end{align*}
Two complex 5-dimensional complete intersections with the same $d,p_1,p_2,e$ are homeomorphic. Furthermore, by Theorem \ref{Traving}, if $\nu_2(d)\geqslant 6.5$ and $\nu_3(d)\geqslant 3.75$, then they will be diffeomorphic.
\begin{example}
Consider the following Table \ref{Diff5},
\begin{table}[h]
\renewcommand{\arraycolsep}{3pt}
\(
\begin{array}{l|*{6}l} \hline
\underline{d} & s_1 & s_2 & s_3 & s_4 & s_5 & d\\ \hline
(112,93,91,50,45,20,18) & 429 & 34723 & 3192813 & 311347699 & 31322739669 & 767763360000\\
(108,105,78,62,35,25,16) & 429 & 34723 & 3192813 & 311347699 & 31322739669 & 767763360000 \\ \hline
(54,48,30,30,13,11,11,4) & 201 & 7447 & 326979 & 15489571 & 763263411 & 14677977600\\
(55,44,39,18,18,16,6,5) & 201 & 7447 & 326979 & 15489571 & 763263411 & 14677977600\\ \hline
\end{array}
\)
\vskip .2cm
\caption{Diffeomorphic 5-dim complete intersections}\label{Diff5}
\end{table}
the two pairs of multidegrees have the same power sums $s_{1},\dots, s_{5}$ and total degree $d$, then it is easily deduced that the corresponding complex 5-dimensional complete intersections have the same Pontrjagin classes and Euler characteristic.
Since total degrees have the following prime factorization,
\begin{align*}
& 767763360000=2^8\times 3^5\times 5^4\times 7^2\times 13\times 31,\\
& 14677977600=2^9\times 3^6\times 5^2\times 11^2\times 13,
\end{align*}
the corresponding two pairs of complete intersections are diffeomorphic.
\end{example}

\begin{example}
The following table gives a diffeomorphic complete intersections with distinct codimensions,
where the diffeomorphism comes from the prime factorization of total degree: $d=104626080000=2^8\times 3^7\times 5^4\times 13\times 23$.
\begin{table}[h]
\(
\begin{array}{l|llll} \hline
X_5(\underline{d}) ~(\textrm{codim=} 9, 10) & d & p_1 & p_2 & e/d \\ \hline
X_5(52,50,30,27,23,18,6,5,4) & 104626080000 & -7748 & 37660770 & -12876778992\\
X_5(54,46,36,25,20,15,13,3,2,2) & 104626080000  & -7748 & 37660770 & -12876778992\\ \hline
\end{array}
\)
\vskip .2cm
\caption{Diffeomorphic 5-dim complete intersections with distinct codimensions}
\end{table}
\end{example}

\subsection{Complex $6$-dimensional complete intersections}
For complex $6$-dimensional complete intersection $X_6(d_1,\dots,d_r)$, the related topological characteristic classes are as follows:
\begin{align*}
 d & =d_1\times\cdots\times  d_r, \\
 p_1 & =7+r-s_2, \\
 p_2 & =\frac{1}{2}{\big [}(7+r-s_2)^2-(7+r-s_4)\big{]},  \\
 p_3 & =\frac{1}{6}\left[(7+r-s_1)^3-3(7+r-s_1)(7+r-s_2)+2(7+r-s_3)\right],  \\
 e & =\frac{d}{6!}\big{[}(7+r-s_1)^6-15(7+r-s_1)^4(7+r-s_2)+40(7+r-s_1)^3(7+r-s_3)\nonumber \\
 & \quad -90(7+r-s_1)^2(7+r-s_4)
 +45(7+r-s_1)^2(7+r-s_2)^2  \\
  & \quad -120(7+r-s_1)(7+r-s_2)(7+r-s_3)+144(7+r-s_1)(7+r-s_5)\\
& \quad -15(7+r-s_2)^3+90(7+r-s_2)(7+r-s_4)\\
& \quad +40(7+r-s_3)^2-120(7+r-s_6) \big ],
\end{align*}
Two complex 6-dimensional complete intersections with the same $d,p_1,p_2,p_3,e$ are homeomorphic. Furthermore, by Theorem \ref{Traving}, if $\nu_2(d)\geqslant 7.5$ and $\nu_3(d)\geqslant 4.25$, then they will be diffeomorphic.

\begin{example}\label{exam6&7}
Consider the following two multidegrees \footnote{The multidegrees in Examples \ref{exam6&7} and \ref{Guo1} are introduced to the author by Guo Xianqiang on mathematical BBS: \url{http://bbs.emath.ac.cn/thread-5853-1-1.html}.}
\begin{align*}
\underline{d} & =(116,114,96,78,59,55,50,40,32,22,13,9),\\
\underline{d}^{\prime} & =(118,110,100,72,64,57,48,39,29,26,11,10),
\end{align*}
They have the same total degree $d$, powers sums $s_1,\dots,s_7$.
\begin{align*}
d & = \fnum{52933656400035840000},\\
s_1 & =\fnum{684},\\
s_2 & =\fnum{54116},\\
s_3 & =\fnum{5008824},\\
s_4 & =\fnum{503305604},\\
s_5 & =\fnum{52970710824},\\
s_6 & =\fnum{5730100991396},\\
s_7 & =\fnum{630552267588024}.
\end{align*}
Then the corresponding two complete intersections have the same total degree, Pontrjagin classes and Euler characteristic, so they are homeomorphic.
Since the total degree has the prime factorization
$d=2^{19}\times 3^{5}\times 5^{4}\times 11^2\times 13^2\times 19\times
29\times 59$,
the corresponding two complete intersections
\begin{align*}
& X_6(116,114,96,78,59,55,50,40,32,22,13,9),\\
& X_6(118,110,100,72,64,57,48,39,29,26,11,10),
\end{align*}
are diffeomorphic.
\end{example}

\subsection{Complex $7$-dimensional complete intersections} \label{subsec7-dim}
For complex $7$-dimensional complete intersection $X_7(d_1,\dots,d_r)$, the related topological characteristic classes are as follows:
\begin{align*}
 p_1 & =8+r-s_2, \\
 p_2 & =\frac{1}{2}{\big [}(8+r-s_2)^2-(8+r-s_4)\big{]},  \\
 p_3 & =\frac{1}{6}\left[(8+r-s_1)^3-3(7+r-s_1)(8+r-s_2)+2(8+r-s_3)\right],  \\
 e & =\frac{d}{7!}\big{[}(8+r-s_1)^7-21(8+r-s_1)^5(8+r-s_2)+70(8+r-s_1)^4(8+r-s_3)\nonumber \\
 & \quad -210(8+r-s_1)^3(8+r-s_4)
 +105(8+r-s_1)^3(8+r-s_2)^2  \\
  & \quad -420(8+r-s_1)^2(8+r-s_2)(8+r-s_3)+504(8+r-s_1)^2(8+r-s_5)\\
& \quad -105(8+r-s_1)(8+r-s_2)^3+630(8+r-s_1)(8+r-s_2)(8+r-s_4)\\
& \quad +280(8+r-s_1)(8+r-s_3)^2-840(8+r-s_1)(8+r-s_6) \\
& \quad +210(8+r-s_2)^2(8+r-s_3)-504(8+r-s_2)(8+r-s_5)\\
& \quad -420(8+r-s_3)(8+r-s_4)+720(8+r-s_{7})\big ],
\end{align*}
Two complex 7-dimensional complete intersections with the same $d,p_1,p_2,p_3,e$ are homeomorphic. Furthermore, by Theorem \ref{Traving}, if $\nu_2(d)\geqslant 8.5$ and $\nu_3(d)\geqslant 4.75$, then they will be diffeomorphic.
\begin{example}\label{exam7&6}
Consider the multidegrees in Example \ref{exam6&7}, it is evidently
that the corresponding two 7-dimensional complete intersections have the same total degree, Pontrjagin classes and Euler characteristic. Thus complete intersections
\begin{align*}
& X_7(116,114,96,78,59,55,50,40,32,22,13,9),\\
& X_7(118,110,100,72,64,57,48,39,29,26,11,10),
\end{align*}
are diffeomorphic by Theorem \ref{Traving}.\hfill\qed
\end{example}

\begin{example}\label{Guo1}
Consider the following two multidegrees 
\begin{align*}
& (596,592,556,520,480,450,438,423,408,404,381,369,327,312,300) \\
& (600,584,564,508,492,447,436,417,416,400,390,360,333,306,303)
\end{align*}
They have the same total degree $d$, powers sums $s_1,\dots,s_{7}$.
\begin{align*}
d & = \fnum{2895548222951602337839765820276736000000},\\
s_1 & =\fnum{6556},\\
s_2 & =\fnum{2994164},\\
s_3 & =\fnum{1424846116},\\
s_4 & =\fnum{703565690996},\\
s_5 & =\fnum{358728296599276},\\
s_6 & =\fnum{187921032698324444},\\
s_7 & =\fnum{100667444609447734036}.
\end{align*}
Furthermore $\nu_2(d)=28$ and $\nu_2(d)=13$,
so the corresponding two complete intersections
\begin{align*}
& X_7(596,592,556,520,480,450,438,423,408,404,381,369,327,312,300), \\
& X_7(600,584,564,508,492,447,436,417,416,400,390,360,333,306,303)
\end{align*}
are diffeomorphic.\hfill\qed
\end{example}

\begin{example}\label{exam11}
For the following two multidegrees \footnote{These two multidegrees were firstly appeared in \url{http://www.emath.ac.cn/florilegium/r011n30.htm} 
by Guo Xianqiang.}
\begin{align*}
\underline{d} & =(12, 16, 22, 26, 28, 45, 58, 59, 65, 69, 81, 85, 86, 91, 105, 106, \\
& \hskip 1cm 108, 128, 132, 134, 144, 156, 168, 192, 200, 214, 242, 250, 272, 274),\\
\underline{d}^{\prime} & =(13, 14, 24, 25, 29, 43, 54, 64, 66, 72, 78, 84, 88, 90, 96, 107, 121, \\
& \hskip 1cm 125, 130, 136, 137, 162, 170, 182, 210, 212, 236, 256, 268, 276),
\end{align*}
They have the same total degree $d$, powers sums $s_1,\dots,s_{11}$.
\begin{align*}
d & = \fnum{43548968602421704369217485857781603627739564212224000000000},\\
s_1 & =\fnum{3568},\\
s_2 & =\fnum{601432},\\
s_3 & =\fnum{120991744},\\
s_4 & =\fnum{26887338904},\\
s_5 & =\fnum{6339294665608},\\
s_6 & =\fnum{1550802794333392},\\
s_7 & =\fnum{388678376991878944},\\
s_8 & =\fnum{99057950894851518184},\\
s_{9} & =\fnum{25552911575591712680248},\\
s_{10} & =\fnum{6651777099183876642569152},\\
s_{11} & =\fnum{1743797004813063555915251344}.
\end{align*}
Then the corresponding two complete intersections are homeomophic.
Since $\nu_2(d)=51$ and $\nu_2(d)=18$,
the corresponding two complete intersections
\begin{align*}
& X_7(12, 16, 22, 26, 28, 45, 58, 59, 65, 69, 81, 85, 86, 91, 105, 106, \\
& \hskip .7cm 108, 128, 132, 134, 144, 156, 168, 192, 200, 214, 242, 250, 272, 274),\\
& X_7(13, 14, 24, 25, 29, 43, 54, 64, 66, 72, 78, 84, 88, 90, 96, 107, 121, \\
& \hskip .7cm 125, 130, 136, 137, 162, 170, 182, 210, 212, 236, 256, 268, 276),
\end{align*}
are diffeomorphic.\hfill\qed
\end{example}

\begin{remark}
Note that, just like Examples \ref{exam6&7} and \ref{exam7&6}, for the multidegrees appeared in Subsection \ref{subsec7-dim}, the corresponding complex 6-dimensional complete intersections are also diffeomorphic.
\end{remark}
\begin{remark}
For higher dimensional complete intersections, e.g. $n=8,9,10,11$, multidegrees in Example \ref{exam11} are still valid to satisfy that the $n$-dimensional complete intersections have the same total degree, Pontrjagin classes and Euler characteristic.  However, because of the unknown classifications of complete intersections for higher dimensions, we can not give more examples of diffeomorphic or homeomorphic complete intersections.
\end{remark}

\noindent{\bf Acknowledgement}\ \ The first author would like to thank Guo Xianqiang for his valuable and helpful discussions.

\end{document}